\documentclass[11pt,reqno]{amsart}
\usepackage{amssymb, amsaddr}
\input amssym.def
\usepackage{amsmath, amsfonts}
\usepackage{amssymb}
\usepackage{amscd, mathtools}
\usepackage[mathscr]{eucal}
\usepackage[usenames]{color}
\usepackage{fancyhdr}
\usepackage{textpos}

\newfont{\cyrr}{wncyr10}

\newcommand{\Q}{{\mathbb Q}}
\newcommand{\R}{{\mathbb R}}
\newcommand{\C}{{\mathbb C}}

\newcommand{\K}{{\mathbf K}}

\renewcommand{\L}{{\mathbf L}}

\newtheorem{thm}{Theorem}

\newtheorem{lem}[thm]{Lemma}
\newtheorem{cor}[thm]{Corollary}

\newtheorem{rmk}[thm]{Remark}

\begin{document}

\title{A note on Dedekind zeta values at $1/2$} 

\author{Neelam Kandhil}

\address{The Institute of Mathematical Sciences, \\
	Homi Bhabha National Institute,\\ 
	Chennai, Tamil Nadu  600 113, India}

\email{neelam@imsc.res.in}


\begin{abstract}   
	For a number field $\K$, let $\zeta_{\K}(s)$ be the Dedekind zeta function associated to $\K$.
	In this note, we study non-vanishing and transcendence of $\zeta_{\K}$ as well
	as its derivative $\zeta_{\K}'$ at $s= 1/2$. En route, we strengthen a result 
	proved by Ram Murty and Tanabe in {\it J. Number Theory, 2016.} 
\end{abstract}
\maketitle

\section{Introduction}

\medskip

 \noindent
For a number field $\K$, let $\zeta_{\K}$ be the Dedekind zeta function associated to $\K$.
The non-vanishing of  $\zeta_{\K}(s)$ at $s=1/2$ is a deep arithmetic question. 
Armitage \cite{Arm} gave examples of number fields  
 $\K$ for which $\zeta_\K(1/2)= 0$. On the other hand, it is believed that $\zeta_{\K}(1/2) \neq 0$  
 when $\K$  is an  $S_n$ - number field, that is, a number field of degree $n$ whose normal closure has Galois group $S_n$ over $\Q$. 
  Furthermore, very little is known about 
 the transcendental nature of the non-zero values of  $\zeta_{\K}(1/2)$. For instance, one has
 $$
 \zeta(1/2) = \frac{1}{ 1 - \sqrt{2}} \sum_{n=1}^{\infty} \frac{(-1)^{n-1}}{\sqrt{n}} 
 \approx -1.46035450880\cdots ,
 $$
 where $\zeta$ is the Riemann zeta function.
 
 \smallskip
 \noindent
 
 In this connection, one has a classical conjecture of Dedekind which asserts that
 if $\L/\K$ is an extension of number fields, then $\zeta_\K(s)$ divides 
 $\zeta_\L(s)$, in other words, the function $\zeta_\L(s) / \zeta_\K(s)$ is entire. This conjecture is open in general,
 but holds when $\L/\K$ is Galois, thanks to the works of Aramata and Brauer. 
 Also the celebrated Artin's conjecture for holomorphicity of his $L$-functions 
 will establish Dedekind's conjecture.

 If $\L/{\Q}$ is a Galois number field with Galois group $S_n$, then
  $\L$ contains a quadratic subfield $\K$. Dedekind's conjecture
 ensures that vanishing of  $\zeta_\K(1/2)$ will ensure vanishing of  $\zeta_\L(1/2)$.
 
In this note, we study  various aspects of the  derivative $\zeta_\K'(s)$ 
at $s= 1/2$. As discussed above, study of these circle of questions for quadratic fields 
merits special attention. We note that for quadratic fields $\K$, the non vanishing of $\zeta_\K(1/2)$ is equivalent to the non vanishing of $L(1/2, \chi)$ for quadratic character $\chi$. Non vanishing of such $L(1/2, \chi)$ has been conjectured by Chowla.
We begin with the following theorem for quadratic fields.

\begin{thm}\label{thm1}  Let  $\K$ and $\L$ be distinct quadratic fields
such that  $\zeta_\K(1/2) \zeta_\L(1/2) \neq 0$. Then
$$\frac{\zeta_\K'(1/2)}{ \zeta_\K(1/2)} \neq  \frac{\zeta_\L'(1/2)}{\zeta_\L(1/2)}.$$
\end{thm}

The above theorem in particular shows that there is at most one quadratic field $\K$  for
which $\zeta_\K(1/2) \ne 0$ while
the derivative $\zeta_\K'(1/2) = 0$. We then prove the following quantitative theorem 
where the existence of the fictitious exception alluded to above is ruled out for both 
quadratic and cubic number fields.

\begin{thm}\label{thm2}  
Let  $\K$ be an algebraic number field with degree $\leq 3$. Then 
$$ 
\zeta_\K(1/2) = 0  \iff  \zeta_\K'(1/2) = 0.
$$
\end{thm}

As  a corollary, we have the following courtesy of the seminal 
work by K. Soundararajan~\cite{Sound}.

\begin{cor}\label{cor3} 
For  at least $87.5 \%$ of quadratic number fields  $\K$ 
with discriminant $8d$ with odd positive square-free $d$, one has
$\zeta_\K'(1/2) \ne 0$.
\end{cor}

Let us very briefly describe the context as well as  content  of the work of Soundararajan indicated above.
The Generalised Riemann Hypothesis (GRH) does not preclude the possibility that $L(1/2, \chi) = 0$ for
some primitive Dirichlet character $\chi$. But it is believed that  there is no rational linear relation between the 
ordinates of the non-trivial zeros of the Dirichlet $\L$-functions and consequently,  $L(1/2, \chi) $ is expected
to be non-zero for any primitive Dirichlet character $\chi$. In particular when $\chi$ is a quadratic character, 
this seems to have been conjectured first by Chowla \cite{chowla} as indicated earlier.  In his outstanding work \cite{Sound}, 
Soundararajan showed that for at least $87.5 \%$ of the odd square-free integers $d \ge 0$, 
$L(1/2 ,\chi_{8d}) \ne 0$. Here for a fundamental discriminant (discriminant of some quadratic number field) $d$,
$\chi_d(n) := \left( \frac{d}{n} \right)$ where $\left( \frac{d}{n} \right)$ denotes 
the Kronecker symbol. Results along  this direction were  obtained 
earlier in \cite{balu}, \cite{isa} and \cite{jut}.

We now have the following theorem for higher degree number fields.

\begin{thm}\label{thm4} Le $\K$ be an algebraic number field of degree $n > 3$ such that the absolute value of its discriminant $|d_\K| \in \R\setminus [(44.763)^{n},(215.333)^{n}]$. Then $$ \zeta_\K(1/2) = 0  \iff  \zeta_\K'(1/2) = 0 . $$
\end{thm}

We now consider the analogous question for Galois number fields.

\begin{thm}\label{thm5}  Consider the following sets.
$$X = \{ \K\,\, {\rm Galois} : \,\, \K \subset \R, \,\, \zeta_\K(1/2) \ne 0,  \,\,\,\, \zeta_\K'(1/2) = 0\}$$
$$Y = \{ \K\,\, {\rm Galois}:  \,\, \K \not\subset \R, \,\, \zeta_\K(1/2) \ne 0,  \,\, \,\, \zeta_\K'(1/2) = 0\}.$$
Then at least one of the sets $X$ and $Y$ is empty.
Furthermore, there are at most finitely many abelian number fields 
for which $\zeta_\K'(1/2) = 0$ but $\zeta_\K(1/2)~\ne~0$.
All such number fields (if exist) have degree less than $46369$.
\end{thm}

\begin{rmk}
Suppose $\zeta_\K(1/2) \ne 0$ and $\zeta_\K'(1/2) = 0,$ then degree of $\K/\Q$ is precisely
 $ \frac{\log |d_\K|}{ \pi/2 + \log 8\pi + \gamma}$
 and $\frac{\log |d_\K|}{\log 8\pi + \gamma}$ in case of totally real and totally complex Galois number fields respectively, where $|d_\K|$ denotes the absolute discriminant of $\K$ and $\gamma$ is the ubiquitous Euler's constant.
\end{rmk}

The above theorem refines a result of Ram Murty and Tanabe \cite [Cor 3.9]{MT}.
Investigations similar to ours for Elliptic curves over $\Q$ as well as Modular forms
were initiated by Gun, Murty and Rath \cite{GMR}.  Furthermore in \cite{MT}, it has been proved that there are 
only finitely many {\it abelian totally real number fields} $\K$ for which $\zeta_\K(1/2) \ne 0$ while
the derivative $\zeta_\K'(1/2) = 0$. One of our objectives in this note was
to  further this line of investigation to arbitrary number fields, obtain some quantitative results 
 and finally study transcendental nature of these 
deeply mysterious numbers. In particular, we use Baker's seminal theorem (see \cite{GMR}, \cite{GMR18}, \cite{MM} and \cite{MP} for some other applications of 
Baker's theorem).
In this context, we have the following theorem.

\begin{thm}\label{thm6}  
Let  $\K$ and $\L$ be distinct algebraic number
fields of degree $n$ and $m$ respectively and  $\zeta_\K(1/2) \zeta_\L(1/2) \neq 0$ .
If one of the two following conditions hold
 \begin{enumerate}
 \item
$$
|d_\K|^m \neq |d_\L|^n;
$$ 
\item
$$
m \frac{\zeta_\K'(1/2)}{\zeta_\K(1/2)} \ne n  \frac{\zeta_\L'(1/2)}{\zeta_\L(1/2)},
$$
\end{enumerate}
then at least one of the following two numbers
$$\frac{\zeta_\K'(1/2)}{ \zeta_\K(1/2)}  \,\,~~ \,{\rm and } \, \,\, ~  \,\, \frac{\zeta_\L'(1/2)}{\zeta_\L(1/2)}$$ 
is transcendental.
\end{thm}

Now, we obtain the following interesting corollaries from Theorem \ref{thm6}.

\begin{cor}\label{cor7}  
Let $n$ be a positive integer. Then the set
$$
\left\{ \frac{\zeta_\K'(1/2)}{ \zeta_\K(1/2)} : \zeta_\K(1/2) \ne 0,  [\K : \Q] = n \right\}
$$
has at most one algebraic number. Furthermore,
$$\frac{\zeta_\K'(1/2)}{ \zeta_\K(1/2)} - \frac{n}{2}(\log 8\pi + \gamma)$$
is a transcendental number.
\end{cor}

Two non-zero integers $u$ and $v$ are said to be multiplicatively independent if for integers $n$ and $m$, $u^n = v^m$ implies $n=m=0$.
In this context, we deduce the following corollary for arbitrary degree number fields.

\begin{cor}\label{cor8}  
Let $\mathcal{F}$ be a family of number fields 
with pairwise multiplicatively independent discriminants. If $\zeta_\K(1/2) \ne 0$
for every $\K \in \mathcal{F}$, then the set
$$
\left\{ \frac{\zeta_\K'(1/2)}{ \zeta_\K(1/2)} : \K \in \mathcal{F} \right\}
$$
has at most one algebraic number.  
\end{cor}

We refer to \cite{MM} and \cite{MP} for investigations related to
the non-vanishing of derivatives of $\zeta_\K(s)$ and $L(s, f)$ at $s=1$, where $f$ is a periodic arithmetic function.

\smallskip

\section{Preliminaries}

\medskip
We fix some useful notations which will be used throughout this paper. For a number field $\K$,

$d_\K$ =  discriminant of $\K$,
 
$r_1$ = number of real embeddings of $\K$,
 
$r_2$ = number of non-conjugate complex embeddings of $\K$,
 
 $n = r_1 + 2r_2$ =  degree of $\K$ over $\Q$,

$\Gamma$ = gamma function,
 
$\gamma$ =  Euler–Mascheroni constant $\approx 0.577\cdots.$ 
\smallskip

We recall some relevant facts about the Dedekind zeta function $\zeta_\K(s)$ associated to a number field $\K$.
$\zeta_\K(s)$ initially given by the following Dirichlet series
$$
\zeta_\K(s) = \sum_{\mathcal{I} \neq 0} \frac{1}{N(\mathcal{I})^s}
$$
for $\Re(s) > 1$ has a meromorphic continuation to the complex plane 
with a simple pole at $s=1$.
Furthermore, the function
$$
Z_\K(s) := \Gamma_{\C}(s)^{r_2} \Gamma_{\R}(s)^{r_1} \zeta_\K(s)
$$
extends meromorphically to the complex plane with simple poles at 
$s=0$ and $s=1$ and satisfies
the functional equation $Z_\K(s) = |d_\K|^{1/2-s} \,\,Z_\K(1-s)$.
Also $\Gamma_{\C} (s) = (2 \pi)^{-s} \Gamma(s),
  \Gamma_{\R} (s) 
= \pi^{-s/2}  \Gamma(s/2)$. But  we shall use the 
following version of the functional equation 
which is amenable for our purpose, namely
\begin{equation*}
\zeta_\K(1-s) = A_\K(s) \zeta_\K(s)
\end{equation*}
for $s \in \mathbb{C}\setminus\{1 \}$ (see \cite[p. 467]{JN}, for instance) with the factor
$$
A_\K(s)
:= 
|d_\K|^{s-1/2} \left(\cos{\frac{\pi s}{2}}\right)^{r_1 + r_2} 
\left( \sin\frac{\pi s }{2} \right)^{r_2} \left( 2(2\pi)^{-s} \Gamma(s) \right)^{n}.
$$ 

We now list some transcendental pre-requisites required for our work.
We shall need Gelfond-Schneider Theorem which states the following.

\begin{thm}\label{GS}
\cite{GE, Nes} 
If $\alpha$ and $\beta$ are non-zero algebraic
numbers with $\beta \ne 1$ and $\log \alpha/ \log \beta \notin \Q,$ then $\log \alpha/ \log \beta $ is transcendental.
\end{thm}

We now record the following application of Baker's  seminal theorem on linear forms in logarithms 
of algebraic numbers \cite[Thm 2.1]{AB}.

\begin{lem}\label{baker}
\cite[p. 154, Lemma 25.4]{MR}  
Let $\alpha_1,\alpha_2, \cdots, \alpha_n $ be positive algebraic numbers. 
If $c_0,c_1, \cdots, c_n$ are algebraic numbers with $c_0 \ne 0$, then 
$$
c_0 \pi + \sum_{j=1}^{n} c_j \log \alpha_j
$$
is a transcendental number.
\end{lem}

Another important lemma required to prove our results
 is the following application of  
Lindemann-Weierstrass theorem \cite{Nes}.

\begin{lem}\label{LW}\cite[Cor 1.3]{Nes}  
If $\alpha$ is an algebraic number different from 
$0$ and $1$, then $\log \alpha$ is a transcendental number where $\log$ 
denotes any branch of logarithmic function.
\end{lem}

Now we quickly recall the discriminant of a quadratic number field $\K$.
Let $d$  be a square-free integer, then the discriminant $d_{\K}$ of the field $\K =\Q({\sqrt {d}})$ is
\begin{equation*}
d_{\K} =
\begin{cases}
d & \text{if $d \equiv $ 1 ~ (mod 4)}\\
 
 4d & \text{if $d \equiv $ 2, 3 ~ (mod 4).}
    \end{cases}       
\end{equation*}

In 1984, Ram Murty proved the following elegant result on  lower bounds of  discriminants of abelian number fields in terms of their degrees. More precisely,

\begin{thm}\label{disc} \cite[Cor 2]{MU} Let $\K$ be an abelian extension of $\Q$ of degree $n$. Then, 
$$
\log |d_\K| \geq \frac{n\log n}{2}.
$$
\end{thm}

We also record the following classical theorem which follows directly from the Minkowski's bound.

\begin{thm}\label{HM} \cite[\S 4.3, Thm 1]{PS} 
For any number field $ \K \ne \mathbb{Q} $, $ |d_\K| > 1$. 
\end{thm}

An essential  ingredient for the proof of  Theorem \ref{thm5} is the following result due to Hermite (see, for instance, \cite[ch. 3]{JN}, \cite{PS}).

\begin{thm}\label{HE} 
There exist only finitely many number fields with bounded discriminant.
\end{thm}

Finally, we state the following deep  theorem of Soundararajan which we shall need  to prove Corollary \ref{cor3}.

\begin{thm}\label{sound} \cite[Thm 1]{Sound}
For at least $87.5\%$ of
the odd square-free integers $d\geq0$,
 $L(1/2,\chi_{8d}) \ne 0$.
\end{thm}

\smallskip

\section{\bf{Proofs of the Main Theorems}}

\medskip

As indicated in the earlier section, we shall work with  the following functional equation of $\zeta_\K(s)$ for $s \in \mathbb{C}\setminus\{1\}$  \cite[p. 467]{JN}
\begin{equation}\label{e-1}
\zeta_\K(1-s) = A_\K(s) \zeta_\K(s)
\end{equation}
with the factor
 $A_\K(s):= |d_\K|^{s-1/2} \left(\cos{\frac{\pi s}{2}}\right)^{r_1 + r_2} \left(\sin\frac{\pi s }{2}\right)^{r_2} \left(2(2\pi)^{-s} \Gamma(s) \right)^{n}.$
 \\
 Let us  begin with an easy, but important  observation  that $A_\K(1/2)$ = 1.
 \\
 We now differentiate \eqref{e-1} w.r.t $s$ and substitute at $ s = 1/2$ to obtain
 \begin{equation}\label{e-2}
\zeta_\K'\left(1/2\right) = -(1/2) A_\K'\left(1/2\right) \zeta_\K\left(1/2\right).
\end{equation}
On the other hand, taking the logarithmic derivative of $A_\K(s)$ we obtain
$$
\frac{A_\K'(s)}{A_\K(s)} = \log |d_\K| - \frac{\pi}{2}(r_1 + r_2) \tan \frac{\pi s}{2} + r_2 \frac{\pi}{2} \cot \frac{\pi s}{2} - n\log 2\pi + n \frac{\Gamma'(s)}{\Gamma(s)},
$$
where $\log $ is the natural logarithm.\\
Since the value of digamma function $ \frac{\Gamma'(s)}{\Gamma(s)}$ at $s= 1/2$ is $ -\gamma - 2\log 2$ (see \cite{SO}, p. 427), we have
\begin{equation}\label{e-3}
A_\K'(1/2) = \log |d_\K| - r_1 \frac{\pi}{2}- n(\log 8\pi + \gamma).
\end{equation}

\noindent
{\bf{3.1.}} {\bf{\textit{ Proof of Theorem \ref{thm1}}}}

\smallskip

\noindent
Let $ \K = \mathbb{Q}(\sqrt{d_1})$ and $ \L = \mathbb{Q}(\sqrt{d_2})$,
where $ d_1$ and $ d_2$
are distinct square-free integers.
Since $\K$ and $\L$ are distinct quadratic fields, we have $ d_\K \ne d_\L$.\\
Using \eqref{e-2} and \eqref{e-3}, we obtain
$$
-2\left(\frac{\zeta_\K'(1/2)}{\zeta_\K(1/2)} -  \frac{\zeta_\L'(1/2)}
{\zeta_\L(1/2)} \right) = 
\log \frac{|d_\K|}{|d_\L|} + \frac{\pi}{2}(r^{(\L)}_1 - r^{(\K)}_1),
$$
where $ r^{(\L)}_1$ and $r^{(\K)}_1$ denote the number of real embeddings of $\L$ and $\K$ respectively.\\
It follows from Theorem \ref{GS} that $e^{\pi}$ is a transcendental number.
So if $r^{(\L)}_1 - r^{(\K)}_1 \ne 0$, then  the right hand side 
of the above equation is non-zero
by Theorem \ref{GS}. On the other hand  if $r^{(\L)}_1 - r^{(\K)}_1 = 0$, 
then the right hand side 
of the above equation is actually transcendental
by Lemma \ref{LW}.\\
Thus,
$$
\frac{\zeta_\K'(1/2)}{ \zeta_\K(1/2)} - \frac{\zeta_\L'(1/2)}{\zeta_\L(1/2)} 
$$
is non-zero.

We note that our proof along with lemma \ref{baker} gives  a stronger assertion, namely
the number $$
\frac{\zeta_\K'(1/2)}{ \zeta_\K(1/2)} - \frac{\zeta_\L'(1/2)}{\zeta_\L(1/2)} 
$$
is actually transcendental. 
\bigskip

\noindent
{\bf{3.2.}} {\bf{\textit{ Proof of Theorem \ref{thm2}}}}

\smallskip
\noindent
By \eqref{e-2}, it is enough to show that $A_\K'\left(1/2\right) \ne 0 $. 
We have
\begin{equation}\label{e-4}
A_\K'\left(1/2\right) = 0 
\iff |d_\K|=
\exp (r_1 \frac{\pi}{2}+ n(\log 8\pi + \gamma)).
\end{equation}
$A_\K'\left(1/2\right) $ is evidently 
 non-zero for $\K =\mathbb{Q}$.
 In fact, we have
 $\zeta'(1/2) = -3.922\cdots$.\\
So we have the following two cases.
\\
{ \it Case (i)}.\,
Assume $\K$ is a quadratic  field.
So $r_1$ could be either 0 or 2.\\
At $r_1 = 0$, we have
$$
2003 < \exp (r_1 \frac{\pi}{2}+ 2(\log 8\pi + \gamma)) < 2004.
$$
At $r_1 = 2$, we have
 $$
 46368 < \exp (r_1 \frac{\pi}{2}+ 2(\log 8\pi + \gamma)) < 46369.
 $$
Since $d_\K$ is always an integer, $A_\K'(1/2)$ can
never be zero in case
of quadratic fields.\\
{ \it Case (ii)}.\, 
Now we consider $\K$ a cubic field.
So $r_1$ is  either 1 or 3.\\
At $r_1 = 1$, we have
$$
431471< \exp (r_1 \frac{\pi}{2}+
3(\log 8\pi + \gamma))< 431472.
$$
At $r_1 = 3$, we have
 $$ 9984558 < \exp (r_1 \frac{\pi}{2}+ 3(\log 8\pi + \gamma))< 9984559.
 $$
 Since $ d_\K$ is always an integer, $A_\K'(1/2)$ can not be zero in case of cubic fields also.

\bigskip

\noindent
{\bf{3.3.}} {\bf{\textit{ Proof of Corollary \ref{cor3}}}}

\smallskip
\noindent
For a quadratic number field $\K$, we have 
$$
\zeta_\K(s) = \zeta (s) L(s,\chi_{d_\K}), \,  \, \, \Re(s) >1,
$$
where $ L(s,\chi_{d_\K}) :=
\sum_{n = 1}^{\infty} \frac{\chi_{d_\K}(n)}{n^s}$. We refer the reader to \cite[Ch. VII]{EH} and \cite{JN}
for further details. We 
recall that this
generalization of Riemann zeta function also has the analytic 
continuation to whole 
complex plane except $ s=1$.
By uniqueness of analytic continuation of complex functions, one 
could get the same identity for $ s \in \mathbb{C}\setminus\{1\}$.\\
Now, we let $ \K= \mathbb{Q}(\sqrt{2d}) $,
where d is a square-free positive odd integer.
It is easy to see that
discriminant of $\K$ is $8d$ (for instance, see \S{5.3}, \cite{PS}). Hence, 
$$
\zeta_\K(1/2) = \zeta (1/2) L(1/2,\chi_{8d}).
$$
Using Theorem \ref{sound}, we have our desired result.

\bigskip

\noindent
{\bf{3.4.}} {\bf{\textit{ Proof of Theorem \ref{thm4}}}}

\smallskip
\noindent
We show that in the given interval of $d_\K$,
$$
A_\K'(1/2) = \log |d_\K| - r_1 \frac{\pi}{2}- n(\log 8\pi + \gamma) \ne 0.
$$
Using hypothesis, we have
$$
A_\K'(1/2) < n(\log (44.763) - \log 8\pi - \gamma) < 0.
$$
Similarly,
$$
A_\K'(1/2) > n(\log (215.333) - \frac{\pi}{2}- \log 8\pi - \gamma) > 0.
$$
So our result follows from \eqref{e-2}.

\bigskip

\noindent
{\bf{3.5.}} {\bf{\textit{ Proof of Theorem \ref{thm5}}}}

\smallskip

\noindent
If possible, let us assume
that there exist Galois number
fields $\K \in X $ 
and $\L \in Y$ of degree $n$ and $m$ respectively. From  \eqref{e-2}, we have
$$
A_\K'(1/2) = A_\L'(1/2)=0.
$$
From \eqref{e-3}, we have
\begin{equation}\label{e-5}
  nA_\K'(1/2) =
  \log |d_\K|^{1/n} - r^{(\K)}_1 \frac{\pi}{2n}- \log 8\pi - \gamma  
\end{equation}
and
\begin{equation}\label{e-6}
mA_\L'(1/2) = \log |d_\L|^{1/m} - r^{(\L)}_1 \frac{\pi}{2m}- \log 8\pi - \gamma,
\end{equation}
where $r^{(\K)}_1$ and $ r^{(\L)}_1$ denote the number of real embeddings of $\K$ and $\L$ respectively.\\
From \eqref{e-5}  and  \eqref{e-6}, we obtain
$$
\log |d_\K|^{1/n} - r^{(\K)}_1 \frac{\pi}{2n}-\log |d_\L|^{1/m} + r^{(\L)}_1 \frac{\pi}{2m} = 0.
$$

Since Galois fields are the normal extensions of their base fields, so there does not exist any complex embedding in real Galois fields. Similarly, there are no real embeddings in non-real Galois fields. Therefore, $r^{(\K)}_1 = n$ and $r^{(\L)}_1 = 0$.
Hence,
$$
\log \frac{|d_\K|^{1/n}}{|d_\L|^{1/m}} - \frac{\pi}{2} = 0, 
$$
which is a contradiction as $ e^{\pi}$ is transcendental by Theorem \ref{GS}. This completes the first part. 

\smallskip

We now proceed with the second part of Theorem \ref{thm5}. From \eqref{e-3}, we have
$$
A_\K'(1/2) = \log |d_\K| - r_1 \frac{\pi}{2}- n(\log 8\pi + \gamma).
$$
Using Theorem \ref{disc}, we obtain
$$
A_\K'(1/2) \geq n(\log (n)/2 - \pi/2 - \log 8\pi - \gamma) > 0, \,\,\, \forall \,  n \geq 46369.
$$
This implies that for all
$ n \geq 46369$, $\zeta_\K'(1/2) = 0$ if and only if $\zeta_\K(1/2) = 0$. \\
Now we aim to prove that the set
$$ S := \{ \K : \, \zeta_\K'(1/2) = 0,\ \zeta_\K(1/2) \ne 0, \,\   n < 46369 \}$$ has finite cardinality.
By  \eqref{e-2}, we see that 
$$ S  \subseteq S':= \{ \K : \, A_\K'(1/2) = 0, \,\   n < 46369 \}\ .$$
So it is enough to show that the set $S'$ has finite cardinality.
By \eqref{e-4}, we have 
$$
A_\K'\left(1/2\right) = 0  \iff |d_\K|= \exp (r_1 \frac{\pi}{2}+ n(\log 8\pi + \gamma)).
$$
Since $n$ and $r_1$ are bounded in the latter set, the discriminant $d_K$
is also bounded. So $S'$ is a set of number fields with bounded discriminant. Hence, we conclude our result by Theorem \ref{HE}.

\bigskip

\noindent
{\bf{3.6.}} {\bf{\textit{ Proof of Theorem \ref{thm6}}}}

\smallskip

\noindent
From \eqref{e-2} and \eqref{e-3}, we obtain

$$
-2\left(m \frac{\zeta_\K'(1/2)}{\zeta_\K(1/2)} - n  \frac{\zeta_\L'(1/2)}
{\zeta_\L(1/2)} \right) =
\log \frac{|d_\K|^m}{|d_\L|^n} + \frac{\pi}{2}(nr^{(\L)}_1 - mr^{(\K)}_1),
$$
where $r^{(\L)}_1$ and 
$r^{(\K)}_1$ denote the number of real embeddings of $\L$ 
and $\K$ respectively.\\
If $nr^{(\L)}_1 -
mr^{(\K)}_1 \ne  0$,
then the right hand side of the above 
 equation is a 
 transcendental number  by Lemma \ref{baker}.
 On the other hand,  if $nr^{(\L)}_1 -
mr^{(\K)}_1 = 0$,
then the right hand side of the above 
 equation is a 
 transcendental number
 by Lemma \ref{LW}.
 So both real numbers
$$\frac{\zeta_\K'(1/2)}
{ \zeta_\K(1/2)}  \,\,\,  {\rm and } \,\,\,  \frac{\zeta_\L'(1/2)}{\zeta_\L(1/2)}$$
can not be algebraic.

\bigskip

\noindent
{\bf{3.7.}} {\bf{\textit{ Proof of Corollary \ref{cor7} and \ref{cor8}}}}

\smallskip

\noindent
If there exist two distinct numbers
from the set given in Corollary \ref{cor7}, then it would be a contradiction to 
Theorem \ref{thm6}.  So the first statement is a direct consequence of Theorem \ref{thm6}.
Now we prove the
second part of Corollary \ref{cor7}. Combining \eqref{e-2} and \eqref{e-3}, we obtain
$$\frac{\zeta_\K'(1/2)}
{ \zeta_\K(1/2)} - \frac{n}{2}(\log 8\pi + \gamma) = r_1 
\frac{\pi}{4} -(1/2)\log |d_\K|.$$
If $ \K= \mathbb{Q}$ ,
then the right hand side of above equation is $\frac{\pi}{4}$. 
By Theorem \ref{HM},  $|d_\K|> 1$ for all $\K$ different from $\mathbb{Q}$. So
$$
\frac{\zeta_\K'(1/2)}{ \zeta_\K(1/2)} - \frac{n}{2}(\log 8\pi + \gamma)
$$  
is a transcendental number by Lemma \ref{baker} and \ref{LW}.

For the proof of  Corollary \ref{cor8}, note that the hypothesis given on discriminants ensures that the first condition 
of Theorem \ref{thm6} is satisfied.  Consequently, it follows  from Theorem \ref{thm6}.

\smallskip

\bigskip

\section{Concluding remarks}

We believe that for any number field $\K$, 
one should  have
$$ \zeta_\K(1/2) \ne 0 \,\, \implies \,\,\zeta_\K'(1/2) \ne 0 . $$

Such results  hold for Elliptic curves over $\Q$ as well as Modular forms \cite{GMR}. The nature of the functional equation in these set ups are amenable
to deduce the above supposition. The classical bounds between degree and discriminant in our context
do not seem to be strong enough to prove the above supposition, at least through our approach.

Furthermore,  if there does exist a number field $\K$ such that  $ \zeta_\K(1/2) \ne 0$
while  $\zeta_\K'(1/2) =0,$ we shall have $\log \pi + \gamma$ being equal to a 
linear form in logarithm of algebraic numbers, an unlikely possibility
from a transcendental perspective since neither $\log \pi$ nor $\gamma$ is expected to be a Baker period, that is, a $\overline{\Q}$ - linear combination of logarithms of algebraic numbers (see \cite{MR} for  details on Baker periods).

\medskip

\bigskip
\noindent
{\bf Acknowledgments.} I would like to express my sincere gratitude 
to Prof. Sanoli Gun for suggesting the problem to me and guiding me 
throughout the duration of the project. It is my immense pleasure to thank 
Prof. Purusottam Rath for useful discussions and ideas.  I thank the referee for a number of suggestions which improved the exposition. I thank 
IMSc also for providing academic facilities. This research is 
supported by a DAE number theory grant.

\newpage

\section*{Addendum: A letter from J.-P. Serre}
\hspace{12.5cm}December 3, 2022

Dear N.Kandhil,

\hspace{0.8cm}About your arXiv text on zeta(1/2).

The behaviour at $s=1/2$ of zeta and $L$ functions is an old question.
I worked briefly on that topic about 50 years ago. Most people seemed to
believe that $L(1/2)$ is nonzero, and they were surprised when I suggested an
explicit counter-example, related to Martinet's work on extensions of $\Q$
with group the quaternion group; that example is the one of Armitage' s
note.
This gave rise to a very interesting theory, mainly developed by
Fröhlich, in which the central point is the sign of the functional
equation. Incidentally, it is strange that you do not mention that sign. In
case it is $+1$, that shows that $L(1/2)= 0$ implies $L'(1/2)=0$.

About $L(1/2)$ being zero or not, the situation seems  to be
the following :

Let $\rho$ be an irreducible complex linear representation of a Galois
extension of  $\Q$  and let $L(s, \rho)$ be the corresponding L function.
Then $L(1/2,\rho)$ is expected to be $\ne 0$, except in the following case :

(*)  $\rho$ fixes a non degenerate antisymmetric bilinear form, and the
constant of the functional equation of  $L(s,\rho)$ is -1.

[ The condition that $\rho$ fixes a non degenerate bilinear form
is equivalent to the character of rho being real valued. When that
bilinear form is symmetric, then the constant of the functional equation
is $+1$ by a theorem of Fröhlich-Queyrut.]

In case (*), one expects that $1/2$ is a simple zero, hence that $L'(1/2,\rho)$
is nonzero.

These expectations are not based on much evidence : merely
that $L$ functions should not have ``unexplained'' common zeroes
(or even $\Q$-related ones). Maybe you could mention them : I am not sure that
they are in the literature.

Best wishes,

J-P. Serre

\end{document}